\documentclass[psamsfonts,12pt]{amsart}

\usepackage[margin=1in]{geometry}  % set the margins to 1in on all sides
\usepackage{graphicx}              % to include figures
\usepackage{amsmath}               % great math stuff
\usepackage{amsfonts}              % for blackboard bold, etc
\usepackage{amsthm}                % better theorem environments
\usepackage{verbatim}              % to comment out large swathes of text with the comment environment
\usepackage{indentfirst}           % to indent the first paragraph of a section
\usepackage{underscore}
\newtheorem{thm}{Theorem}[section]
\newtheorem{lem}[thm]{Lemma}
\newtheorem{prop}[thm]{Proposition}
\newtheorem{cor}[thm]{Corollary}

\theoremstyle{remark}
       \newtheorem*{rmk}{Remark}
\theoremstyle{remark}

\newcommand{\CC}{\mathbb{C}}
\newcommand{\ZZ}{\mathbb{Z}}
\newcommand{\Mod}[1]{\ (\textup{mod}\ #1)}     
\newcommand{\QQ}{\mathbb{Q}}
\newcommand{\Qbar}{\overline{\QQ}}
\newcommand{\Pone}{\mathbb{P}^1} 
\newcommand{\PoneQbar}{\Pone(\Qbar)}
\newcommand{\PoneQD}{\Pone(\Qbar,d)} 
\newcommand{\PDB}{\Pone(\overline{\QQ},d,B)}
\newcommand{\OS}{\mathcal{O}_S}
 \newcommand{\Orb}{\operatorname{Orb}}

\title[Integral Points of Bounded Degree]{Integral Points of Bounded Degree on the Projective Line and in Dynamical Orbits}
\date{\today}
\author[Joseph Gunther]{Joseph Gunther}
\address{Department of Mathematics, The Graduate Center, City University of New York (CUNY); 365 Fifth Avenue, New York, NY 10016 USA}
\email{JGunther@gradcenter.cuny.edu}
\author[Wade Hindes]{Wade Hindes}
\address{Department of Mathematics, The Graduate Center, City University of New York (CUNY); 365 Fifth Avenue, New York, NY 10016 USA}
\email{whindes@gc.cuny.edu}
\subjclass[2000]{}
\begin{document}

\maketitle
\renewcommand{\thefootnote}{}
\footnote{The first author was partially supported by National Science Foundation grant DMS-1301690.}
\begin{abstract} Let  $D$ be a non-empty effective divisor on $\mathbb{P}^1$. We show that when ordered by height, any set of $(D,S)$-integral points on $\mathbb{P}^1$ of bounded degree has relative density zero.  We then apply this to arithmetic dynamics: let $\varphi(z)\in \Qbar(z)$ be a rational function of degree at least two whose second iterate $\varphi^2(z)$ is not a polynomial. We show that as we vary over points $P\in\mathbb{P}^1(\Qbar)$ of bounded degree, the number of algebraic integers in the forward orbit of $P$ is absolutely bounded and zero on average.
\end{abstract} 

\section{Introduction}
Let $K/\mathbb{Q}$ be a number field and let $S$ be a finite set of places of $K$ containing the archimedean ones. Siegel's theorem is a fundamental result in the study of integral points on curves (we use below the modern language of $(D,S)$-integral points; see Section 2 for more details):

\begin{thm}[Siegel]

Let $C$ be a curve defined over $K$ and let $D$ be a non-empty effective divisor on $C$, also defined over $K$.  Then if $D$ contains at least 3 distinct (geometric) points, any set of $(D,S)$-integral points in $C(K)$ is finite. 
\end{thm}

More classically, one can state this (and later results of this introduction) in terms of an affine embedding: if $C_{/K}\subseteq \mathbb{A}^n$ is an affine curve with at least three points at infinity, then $C$ has only finitely many $K$-points whose coordinates are $S$-integers.

In fact, a curve $C$ can have infinitely many integral $K$-points only if $D$ comprises one or two points and $C$ has genus zero.  Even in this infinite case, integral points are still known to be rare, in the following sense: when we order the $K$-points of $C$ by a Weil height $H(\cdot)$ on $C$, any set of integral points has density zero within the rational points of $C$.

\begin{thm} \label{integralsparse}
Let $D$ be a non-empty effective divisor on $\Pone$ defined over $K$, and let $\mathcal{R}$ be any set of $(D,S)$-integral points in $\Pone(K)$. Then $${\displaystyle \lim_{B \rightarrow \infty}}\; \frac{\#\{P \in \mathcal{R} \ | \ H(P) \leq B\}}{\#\{P \in \Pone(K) \ | \ H(P) \leq B\}} = 0.$$
\end{thm}

See \cite[Chapter 9]{serre} for a proof in the case where $S$ is exactly the set of archimedean places; it also shows the proportion of integral points decreases to zero relatively fast.

More recently, Siegel's theorem has been generalized, beyond integral points defined over a fixed number field, to integral points defined over varying number fields of bounded degree.  This deep finiteness result \cite[Corollary 14.14]{levindegreed} follows from work of Vojta \cite{vojtapaper} and of Song and Tucker \cite{songtucker}:

\begin{thm}[Levin] \label{levinfiniteness}
Let $d \geq 1$ be an integer, let $C$ be a curve defined over $K$, and let $D$ be a non-empty effective divisor on $C$, also defined over $K$.  Then if $D$ contains at least $2d + 1$ points, any set of $(D,S)$-integral points contained in $\{P \in C(\Qbar) \ | \ [K(P):K] \leq d\}$ is finite. 
\end{thm}
\begin{rmk}
Note that when $d=1$, this recovers Siegel's theorem.  See \cite{levinintegral} for an elegant converse, which shows that in one sense, integral points of bounded degree on curves behave better than rational points.  Some special cases of Theorem \ref{levinfiniteness} were known earlier: if $C=\Pone$, it follows from the Thue-Siegel-Roth-Wirsing theorem \cite{wirsing} on Diophantine approximation.  For arbitrary curves, Corvaja and Zannier proved the theorem \cite[Corollary 1]{cz} in the case $d=2$.
\end{rmk}

Given this finiteness result, it's reasonable to ask if a higher-degree analogue of Theorem \ref{integralsparse} holds. In other words: in the cases where $C$ has infinitely many integral points of degree $d$, are they still density zero within the rational points of degree $d$?  Our first theorem shows this is indeed true if $C$ has genus zero and the base field is $\QQ$.  For $d \geq 1$, let $\PoneQD = \{P \in \PoneQbar \ | \ [\QQ(P):\QQ] \leq d \}$.

\begin{thm} \label{densityzero}
Let $D$ be a non-empty effective divisor on $\Pone$ defined over $K$, and let $\mathcal{R} \subset \PoneQD$ be a set of $(D,S)$-integral points.  Then $${\displaystyle \lim_{B \rightarrow \infty}}\; \frac{\#\{P \in \mathcal{R} \ | \ H(P) \leq B\}}{\#\{P \in \PoneQD \ | \ H(P) \leq B\}} = 0.$$ 
\end{thm}

\begin{rmk}
If $d=1$ this is a special case of Theorem \ref{integralsparse}, while if $D$ contains a $\QQ$-point, it follows from \cite[Theorem 1.2]{barroero} or \cite[Theorem 3.5.6]{integralbatman}.  However, in its full generality, Theorem \ref{densityzero} appears to be new.
\end{rmk}

While Theorem \ref{densityzero} is of independent interest, our original motivation for proving it was an application to arithmetic dynamics, which we'll now explain.

In \cite{Silv-Int}, Silverman established the following dynamical corollary to Siegel's theorem (here $\varphi^n = \varphi \circ \dots \circ \varphi$ denotes the $n$th iterate of $\varphi$):  
\begin{thm}[Silverman]{\label{thm:Silverman}} Let $\varphi(z)\in K(z)$ be a rational function of degree at least two. If $\varphi^2(z)$ is not a polynomial, then the forward orbit of $P\in\mathbb{P}^1(K)$,   
\[\Orb_\varphi(P):=\{P,\varphi(P),\varphi^2(P),\dots\},\]
contains only finitely many $S$-integral points. 
\end{thm} 
In light of Theorem \ref{thm:Silverman}, it is natural to ask if the number of integral points in an orbit of $\varphi$ can be uniformly bounded.  We show that the answer is yes, and in fact we strengthen the statement in two ways: first, we allow $\varphi$ to be defined with arbitrary $\Qbar$-coefficients, and second, our bound depends only on the degree of the number field.

Some notation: for $S$ a finite set of places of $\QQ$ containing $\infty$, let $\OS$ denote the integral closure of $\ZZ_S$ within $\Qbar$.  When $S=\{\infty\}$, this is simply the ring of all algebraic integers; the reader is welcome to keep this intuitive example in mind throughout the paper.  Whenever we write $P \in \OS$ for a point of $\PoneQbar$, that means $P$ is of the form $[\alpha:1]$, for $\alpha \in \OS$.
\begin{thm} \label{maxnumberpoints}
Let $\varphi(z) \in \Qbar(z)$ be a rational function of degree at least two and let $S$ be a finite set of places of $\QQ$ containing the archimedean one. Then if $\varphi^2(z)$ is not a polynomial, there exists a constant $N=N(\varphi, d, S)$ such that for any point $P \in \PoneQD$, we have $\#(\Orb_{\varphi}(P) \cap \OS) \leq N$.
\end{thm}
   
Although it's nice to have an upper bound on the number of integral points of $\Orb_\varphi(P)$, one expects that most orbits contain no integers at all. To test this intuition, we study the average number of integral points in orbits as we vary over $P\in\PoneQbar$ of degree at most $d$ and height at most $B$.  In particular, we show that this average tends to zero as the height grows. Moreover, since the choice of $d$ is arbitrary, the following result can be roughly interpreted as: ``a random algebraic number has no integral points in its orbit."

For $P \in \PoneQbar$, let $H(P)$ denote its absolute multiplicative Weil height.  Let $\PDB = \{P \in \PoneQbar \ | \ [\QQ(P):\QQ] \leq d \text{ and } H(P) \leq B\}$; by Northcott's theorem, this is a finite set.  Then we have:

\begin{thm}{\label{average-int-orbit}}
Let $\varphi(z) \in \Qbar(z)$ be a rational function of degree at least two and let $S$ be a finite set of places of $\QQ$ containing the archimedean one. Then if $\varphi^2(z)$ is not a polynomial, $${\displaystyle \lim_{B \rightarrow \infty}} \frac{\sum_{P \in \PDB} \#(\Orb_{\varphi}(P) \cap \OS)}{\#\PDB} = 0$$ for all $d\geq1$. 
\end{thm}
Aside from standard properties of dynamical heights, the two main ingredients in the proof are the finiteness result in Theorem \ref{levinfiniteness} and the density result in Theorem \ref{densityzero}.

The dynamical results in this paper generalize \cite[Theorem 1.2]{wadepaper}, where the second author studies averages over a fixed number field.
\section{Notation and Previous Results}

In this section, we establish some notation for the rest of the paper, and state known results that we'll use.

There are various roughly equivalent definitions of integral points on a variety, so we'll fix a convenient one for our purposes.  For a number field $K$, a variety $X$ defined over $K$, a divisor $D$ on $X$ defined over $K$, and a finite set $S$ of places of $K$ containing the archimedean ones, we say that $\mathcal{R} \subset (X\backslash D)(\Qbar)$ is a set of $(D, S)$-integral points if there exists a global Weil function $\lambda_D$ such that for all places $v$ of $K$ not in $S$, we have $\lambda_{D,v}(P) \leq 0$ for all $P \in \mathcal{R}$.  We refer the reader to \cite[Chapter 1]{vojtabook} and \cite[Chapter 10]{lang} for the definition of a global Weil function as a suitable collection of local Weil functions; we'll simply cite the properties relevant to us as needed.  

The following ``clearing denominators" lemma \cite[Lemma 1.4.6]{vojtabook} will be of particular use:

\begin{lem} \label{cleardenominators}
Let $\mathcal{R}$ be a $(D,S)$-integral set of points on $X$, and let $f$ be a rational function with no poles outside of $D$.  Then there is some constant $b \in K^\times$ such that $b\cdot f(P)$ is integral for all $P \in \mathcal{R}$.
\end{lem}

With the exception of Theorem \ref{densityzero}, however, the above general notion of integrality is only used within the proofs, not in stating our dynamical results.  For that, we need only the simpler definition of $\OS$ given above Theorem \ref{maxnumberpoints}.

We'll use two different notions of height.  First, for a point $P \in \PoneQbar$, let $H(P)$ denote the absolute multiplicative Weil height of $P$, for which there are various equivalent definitions. One of them: let $H([1:0]) = 1$, while for $\alpha \in \Qbar$, let $H([\alpha:1]) = \left(a_d\prod_{\text{conjugates }\alpha' \text{ of }\alpha} \operatorname{max}\{1, |\alpha'|\}\right)^{\frac{1}{[\QQ(\alpha):\QQ]}}$, where $a_d$ is the leading coefficient of $\alpha$'s {\em minimal polynomial over $\ZZ$}; here we obtain the minimal polynomial over $\mathbb{Z}$ from the minimal polynomial over $\QQ$ by multiplying by the smallest positive integer that clears denominators.

Second, to a non-constant rational map $\varphi\in \Qbar(z)$ of degree $r \geq 2$, we can associate a different height function $\hat{h}_{\varphi}$ on $\PoneQbar$ called the (logarithmic) canonical height; see, for instance \cite[\S3.4]{silvermanbook}. The canonical height behaves well under iteration: $\hat{h}_{\varphi}(\varphi^n(P)) = r^n \hat{h}_{\varphi}(P)$.  Moreover, a point $P$ has canonical height $0$ if and only if it is a {\em preperiodic} point for $\varphi$, i.e. its forward orbit is finite.

Both of these heights satisfy the Northcott property \cite[Theorems 3.7 and 3.12]{silvermanbook}: for fixed values of $d$ and $B$, there are only finitely many points of $\PoneQbar$ of degree at most $d$ and height at most $B$.

Next, for $d, B \geq 1$, we defined in the introduction $\PoneQD = \{P \in \PoneQbar \ | \ [\QQ(P):\QQ] \leq d \}$, as well as the (finite) subset $\PDB = \{P \in \PoneQbar \ | \ [\QQ(P):\QQ] \leq d \text{ and } H(P) \leq B\}$. While an easy crude bound on the size of $\PDB$ would suffice for our purposes, it's clarifying to state some recent stronger results.  For $d$ fixed and $B$ increasing, Masser and Vaaler \cite{masservaaler} determined its asymptotic size:

\begin{thm}[Masser-Vaaler]\label{massvaal}
As $B$ grows, the number of elements $\alpha \in \Qbar$ such that $[\QQ(\alpha):\QQ]=d$ and $H(\alpha)\leq B$ is asymptotic to $c_{\QQ,d} \cdot B^{d(d+1)}$, for $c_{\QQ,d}$ an explicit positive constant.
\end{thm}

More generally, they established an asymptotic for points of degree $d$ over arbitrary number fields \cite{masservaaler2}, as well as a power-saving error term, but we won't need that here.  If we restrict attention to $S$-integral points, Barroero \cite[Theorem 1.2]{barroero} showed:

\begin{thm}[Barroero]\label{barr}
Let $S$ be a finite set of places of $\QQ$ containing the archimedean one.  Then as $B$ grows, the number of elements $\alpha \in \OS$ such that $[\QQ(\alpha):\QQ]=d$ and $H(\alpha) \leq B$ is asymptotic to $a_{\QQ, d, S} \cdot B^{d^2}(\operatorname{log }B)^{|S| - 1}$, for $a_{\QQ, d, S}$ an explicit positive constant.
\end{thm}

Barroero also obtained an error term, as well as the asymptotic over an arbitrary base number field rather than just over $\QQ$, but again we won't need that full generality.
\section{Finiteness of integral images}

\begin{lem} \label{numberofpoles}
Let $\varphi(z) \in \CC(z)$ be a rational map of degree $r \geq 2$.  If $\varphi^2(z)$ is not a polynomial, then the number of distinct poles of $\varphi^n(z)$ goes to $\infty$ as $n \rightarrow \infty$. 
\end{lem}

\begin{proof}
If $\infty$ is not a periodic point for $\varphi$, then $\varphi^n$ has at least $r^{n-2}$ poles; see \cite[Exercise 3.37(a)]{silvermanbook}. On the other hand, if $\infty$ is a periodic point of exact order $m$, then note that by our assumption on $\varphi^2$ and \cite[Theorem 1.7]{silvermanbook}, $\infty$ is a fixed point of $\varphi^m$ but not totally ramified.  Hence $\varphi^{mn}$ has at least $r^{n-2}+2$ poles by \cite[Exercise 3.37(b)]{silvermanbook}.  In either case, the statement follows.
\end{proof}

We'll write $x_0$ and $x_1$ for the natural coordinates on the projective line.

\begin{prop} \label{finiteintegral}
Let $\varphi(z) \in \Qbar(z)$ be a rational map with at least $2d+1$ distinct poles and let $S$ be a finite set of places of $\QQ$ containing the archimedean one.  Then there are only finitely many points $P \in \PoneQD$ such that $\varphi(P) \in \OS$.
\end{prop}

\begin{proof}
Let $\mathcal{R}=\varphi(\PoneQD)\cap \OS$, let $K$ be the smallest number field over which the coefficients of $\varphi$ can be defined, and let $S'$ be the set of places of $K$ lying over $S$.  By \cite[(1.3.5)]{vojtabook}, setting $\lambda_{\{\infty\},v}=\frac{1}{[K:\QQ]}\text{log max}(1, || \frac{x_0}{x_1} ||_v)$ for $v \not \in S'$ gives a global Weil function $\lambda_{\{\infty\}}$, which thus shows $\mathcal{R}$ to be a set of $(\{\infty\}, S')$-integral points.  By \cite[Lemma 1.3.3(d)]{vojtabook}, $\lambda_{\{\infty\}} \circ \varphi$ is a global Weil function for the divisor $\varphi^*\{\infty\}$ on $\Pone$.  From our earlier definition of integral points in Section 2, the subset of $\PoneQD$ as in the proposition is immediately seen to be $(\varphi^*\{\infty\}, S')$-integral.  By Theorem \ref{levinfiniteness}, there can only be finitely many such points of degree at most $d$, since $\varphi^*\{\infty\}$ contains at least $2d+1$ distinct points by assumption.
\end{proof}

\section{Upper bounds for orbits}

\begin{prop} \label{maxiterate}
Let $\varphi(z) \in \Qbar(z)$ be a rational map of degree $r \geq 2$ such that $\varphi^2(z)$ is not a polynomial, and let $S$ be a finite set of places of $\QQ$ containing the archimedean one.  Then there exists a constant $N'=N'(\varphi, d, S)$ such that for any non-preperiodic point $Q \in \PoneQD$, we have that $\varphi^n(Q) \in \OS$ implies $n \leq N'$.
\end{prop}

\begin{proof}
By Lemma \ref{numberofpoles}, there exists $N_1$ such that $\varphi^{N_1}$ has at least $2d+1$ distinct poles.  By Proposition \ref{finiteintegral}, there are only finitely many points $P \in \PoneQD$ such that $\varphi^{N_1}(P) \in \OS$.  Let $C$ be the maximum of the canonical heights $\hat{h}_{\varphi}(P)$ of these points.  If $\varphi^n(Q) \in \OS$ for $Q$ as in the proposition and $n \geq N_1$, then $\varphi^{n-N_1}(Q)$ has height at most $C$.  Thus 
\[r^n\hat{h}_{\varphi}(Q) = \hat{h}_{\varphi}(\varphi^n(Q)) = \hat{h}_{\varphi}(\varphi^{N_1}(\varphi^{n-N_1}(Q))) = r^{N_1}\hat{h}_{\varphi}(\varphi^{n-N_1}(Q)) \leq r^{N_1}C.\]  So if $M=M(\varphi,d)$ is the minimal positive value of the canonical height $\hat{h}_{\varphi}$ on $\PoneQD$ (which immediately exists by the Northcott property), then we must have $r^n \leq \frac{r^{N_1}C}{M}$, which proves the proposition.
\end{proof}
\begin{proof}[(Proof of Theorem \ref{maxnumberpoints})]
Proposition \ref{maxiterate} immediately reduces this to finding a uniform bound just for such points $P$ that are preperiodic.  But those points have canonical height 0, so by Northcott there are finitely many of them, since their degree is bounded by assumption.  Thus any bound bigger than the orbit lengths of all the preperiodic $P$'s will suffice.
\end{proof}
\begin{rmk} If $\varphi^2(z)\in \Qbar[z]$ and $\varphi(z)$ is not itself a polynomial, then after a change of variables $\varphi(z)$ has the form $1/z^r$; see \cite[Proposition 1.1]{Silv-Int}. 
\end{rmk}
\section{Average number of integral points in orbits}

\begin{proof}[(Proof of Theorem \ref{densityzero})]
Recall that as $B$ grows, $\#\Pone(\Qbar,d,B)$ is asymptotic to $c_{\QQ,d} \cdot B^{d(d+1)}$ by Theorem \ref{massvaal}.

First, suppose that $D$ contains $\infty=[1:0]$.  Then $\frac{x_0}{x_1}$ is a regular function on the complement of $D$, so by Lemma \ref{cleardenominators}, there is a constant $b \in \Qbar^\times$ such that $\frac{bx_0}{x_1}(P) \in \mathcal{O}_{S'}$ for all $P \in \mathcal{R}$, where $S'$ is the finite set of places of $\QQ$ that the places of $S$ lie over. If we expand $S'$ to a possibly larger finite set of places $T$ of $\QQ$ such that it contains all the places above which $b$ has absolute value less than 1, we see that $\mathcal{R} \subset \mathcal{O}_{T}$.  Thus $\#\mathcal{R} \cap \PDB$ is bounded above by $\# \mathcal{O}_T \cap \PDB$.  Since $\# \mathcal{O}_T \cap \PDB$ is asymptotic to $a_{\QQ,d,T} \cdot B^{d^2}(\text{log }B)^{|T|-1}$ as $B$ grows, the limit in the theorem statement is zero.

Now suppose instead that $D$ doesn't contain $\infty$.  Then it has a finite point $[\beta:1]$.  Consider the rational function $\frac{x_1}{x_0-\beta x_1}$; this is a regular function on the complement of $D$, so by Lemma \ref{cleardenominators}, there is again a constant $b \in \Qbar^\times$ such that $\frac{bx_1}{x_0-\beta x_1}(P)\in \mathcal{O}_{S'}$ for all $P \in \mathcal{R}$. In particular, for all non-infinite $[\alpha:1]\in\mathcal{R}$, we see that $\frac{1}{\alpha-\beta}\in\mathcal{O}_T$, where $T$ is the finite set of places containing $S'$ and all places above which $b$ can have absolute value less than 1. Hence, to prove Proposition \ref{densityzero} it suffices to show that the set of points $\alpha\in\Qbar$ satisfying $\frac{1}{\alpha-\beta}\in\mathcal{O}_T$ and $[\mathbb{Q}(\alpha):\mathbb{Q}]\leq d$ has relative density zero inside $\PoneQD$.  In the special case where one can choose $\beta$ to lie in $\QQ$, we again get an asymptotic bound of a constant times $B^{d^2}(\text{log }B)^{|T|-1}$ from Theorem \ref{barr}.  But for the general case, this doesn't work; instead we'll show relative density zero by sieving out a family of local conditions.

First we make a couple reductions.  Since $\#\Pone(\overline{\QQ},d-1,B)$ is asymptotic to $c_{\QQ,d-1}\cdot B^{d(d-1)}$, one need only handle the set where $[\QQ(\alpha):\QQ]= d$.  Furthermore, by an appropriate form of effective Hilbert irreducibility \cite[Theorem 2.1]{cohen} over $\QQ(\beta)$, we may restrict attention to the subset of $\alpha$'s where $[\QQ(\alpha,\beta):\QQ(\beta)]= [\QQ(\alpha):\QQ]= d$.

To outline the argument going forward, let $\mathcal{P}$ be the set of rational primes that split completely in the Galois closure of $\mathbb{Q}(\beta)$. By the Chebotarev density theorem, $\mathcal{P}$ has positive Dirichlet density. We may discard finitely many primes of $\mathcal{P}$ and assume that $\beta$ is integral at all $p\in\mathcal{P}$ and that $\mathcal{P}$ does not meet $T$. Next, for each $p\in\mathcal{P}$, we can fix a prime ideal $\mathfrak{p}$ of $\mathbb{Q}(\beta)$ lying over $p$, and since there is no residue field extension we may choose an integer $0\leq r_p\leq p-1$ such that $| r_p-\beta |_\mathfrak{p}<1$. 
Now, consider the set
\begin{equation*}{\label{congruence}}
\mathcal{I}_p:=\big\{\alpha\in\Qbar \ \big| \ [\mathbb{Q}(\alpha):\mathbb{Q}]= d, \; | \alpha - r_p | < 1 \text{ for some absolute value }| \cdot |\text{ of }\Qbar \text{ lying over } | \cdot |_\mathfrak{p}\big\}.
\end{equation*}  
For an element $\alpha$ of $\mathcal{I}_p$, we have $|\frac{1}{\alpha-\beta}|=|\frac{1}{(\alpha-r_p)+(r_p-\beta)}| > 1$, and therefore $\frac{1}{\alpha-\beta} \not \in \mathcal{O}_T$.  Thus $\mathcal{R}\cap\mathcal{I}_p = \emptyset$, and so we'll be done if we show that the complement of $\cup_{p \in \mathcal{P}}\mathcal{I}_p$ has density zero. Rather than counting elements of the complement directly, we'll instead bound its size from above, by counting polynomials whose roots lie in the complement.
 
To make this precise, let $f(x)=a_dx^d+a_{d-1}x^{d-1}+\dots+a_0$ be a primitive irreducible polynomial with integer coefficients and write $f(x)=a_d\prod(x-\alpha_i)$ in $\mathbb{C}[x]$. Then we define $\hat{H}(f):=\max\{|a_i|\}$ to be the naive height of $f$. For $\alpha$ a root of $f$, the height of $\alpha$ and the naive height of $f$ are known to be comparable in the following sense: 
\begin{equation}{\label{inequality}} \frac{1}{\sqrt{d+1}}\,H(\alpha)^d\leq \hat{H}(f)\leq {d\choose{\lfloor d/2\rfloor}} H(\alpha)^d;
\end{equation}  
see, for instance \cite[Lemma 1.6.7]{bombierigubler}. Now let $p_1, \dots, p_k$ be the first $k$ primes in $\mathcal{P}$.  Let $\operatorname{Pol}^+(d,B)$ be the set of integer polynomials of degree $d$, with positive leading coefficient, and naive height at most $B$. Lastly, for any $m | p_1\dotsm p_k$, let  
\begin{equation}\label{polyswithcongruence}
\begin{split}
\mathcal{F}_m(d,B):=\#\big\{ f\in \operatorname{Pol}^+(d,B) \ \big| \ \text{for each }p | m, \ a_d\not\equiv0\Mod{p}, f(r_p)\equiv0\Mod{p}\big\}.
\end{split}
\end{equation} 
By counting the number of integers in a residue class in a box (at each prime we're excluding a single congruence class for the leading coefficient, and enforcing one condition on the remaining terms), we see that as $B$ grows,
\begin{equation*} \label{naivecount}
\mathcal{F}_m(d,B)=\prod_{p|m}\Big(\frac{p-1}{p^2}\Big)2^dB^{d+1}+O(B^d)
\end{equation*} 
On the other hand, the irreducible polynomials counted by $\mathcal{F}_m(d,B)$ give rise to elements of $\cap_{p|m}\mathcal{I}_p$; to see this, suppose that $\alpha\in\Qbar$ is a root of an {\em irreducible} integer polynomial $f(x)$ of this type. Then considering the equation
\begin{align*}
0 =a_d\alpha^d+a_{d-1}\alpha^{d-1}+\dots a_0 & =a_d(\alpha-r_p+r_p)^d+a_{d-1}(\alpha-r_p+r_p)^{d-1}+\dots a_0\\
& =a_d(\alpha-r_p)^d+ (da_dr_p+a_{d-1})(\alpha-r_p)^{d-1} + \dots + f(r_p),
\end{align*}
 we see, dividing by $a_d$ to obtain the minimal polynomial equation for $\alpha - r_p$ over $\QQ$, that the norm $N_{\QQ(\alpha)/\QQ}(\alpha-r_p) = \frac{(-1)^df(r_p)}{a_d}$.  Let $| \cdot |$ be an absolute value of $\Qbar$ lying over $| \cdot |_\mathfrak{p}$.  Since the norm, which is the product of the Galois conjugates of $\alpha-r_p$, has absolute value less than 1 under $| \cdot |$, we must have, for some automorphism $\sigma \in \operatorname{Gal}(\Qbar/\QQ)$, that $\sigma(\alpha-r_p)=\sigma(\alpha)-r_p$ satisfies $|\sigma(\alpha)-r_p|<1$, and thus $\sigma(\alpha) \in \mathcal{I}_p$.  Now, by our earlier reduction to the case $[\QQ(\alpha,\beta):\QQ(\beta)]= [\QQ(\alpha):\QQ]= d$, we have that $\operatorname{Gal}(\Qbar/\QQ(\beta))$ acts transitively on the conjugates of $\alpha$, so in fact we can choose $\sigma$ to lie in $\operatorname{Gal}(\Qbar/\QQ(\beta))$.  Thus the absolute value $| \cdot |'$ of $\Qbar$ given by $| \cdot |' = | \cdot | \circ \sigma$ also lies over $| \cdot |_\mathfrak{p}$, and we have $| \alpha - r_p |' < 1$.  So we in fact have $\alpha \in \mathcal{I}_p$.\

Next, let $\mathcal{G}_k(d,B)$ be the number of polynomials in $\operatorname{Pol}^+(d,B)$ not contained in any of the subsets defining $\mathcal{F}_{p_i}(d,B)$ for $i=1,\dots,k$.  By inclusion-exclusion, we have $$\mathcal{G}_k(d,B) = \sum_{m|p_1\dotsm p_k} \mu(m)\mathcal{F}_m(d,B)$$ $$= \prod_{i=1}^k\Bigg(1-\Big(\frac{p_i-1}{p_i^2}\Big)\Bigg)2^dB^{d+1} + O\big(B^d\big).$$

Now we count algebraic numbers.  Let $\alpha \in \Qbar$ have degree $d$, and let $f$ be its minimal polynomial over $\ZZ$.  If $\alpha$ has height at most $B$, then by (\ref{inequality}), we see $f$ has naive height $\hat{H}(f) \leq {d\choose{\lfloor d/2\rfloor}} B^d$.  Thus the number of such $\alpha$ not contained in $\cup_{p \in \mathcal{P}}\mathcal{I}_p$ is bounded above by $$d \cdot \mathcal{G}_k(d,{d\choose{\lfloor d/2\rfloor}}B^d)$$ $$= \prod_{i=1}^k\Bigg(1-\Big(\frac{p_i-1}{p_i^2}\Big)\Bigg)d2^d{d\choose{\lfloor d/2\rfloor}}^{d+1}B^{d(d+1)} + O\big(B^{d^2}\big).$$

We see that the relative density of $\mathcal{R}$ is at most $\frac{1}{c_{\QQ,d}}$ times the constant on $B^{d(d+1)}$ above.  However, now we can let $k$ grow.  Since $\mathcal{P}$ had positive density in the primes, the product above converges to 0 as $k \to \infty$: recall that an infinite product $\prod (1 - a_i)$, with $0 \leq a_i < 1$, converges to 0 if and only if $\sum a_i$ diverges.  In our case, $a_i$ is on the order of $\frac{1}{p_i}$; since the sum of the reciprocals of the primes diverges \cite[Theorema 19]{euler}, the infinite product converges to $0$.\end{proof} 
\begin{cor} \label{imagedensity}
If $\varphi(z) \in \Qbar(z)$ is a non-constant rational function, then $${\displaystyle \lim_{B \rightarrow \infty}} \frac{\#\{P \in \PDB \ | \ \varphi(P) \in \OS\}}{\#\PDB} = 0.$$
\end{cor}

\begin{proof}
This follows directly from Theorem \ref{densityzero}: as noted in the proof of Proposition \ref{finiteintegral}, we have that $\{P \in \PoneQD \ | \ \varphi(P) \in \OS\}$ is a $(\varphi^*\{\infty\},S')$-integral set of points.
\end{proof}  

\begin{proof}[(Proof of Theorem \ref{average-int-orbit})]
By Proposition \ref{maxiterate}, after discarding the finitely many preperiodic points, the numerator is at most $$\sum_{n=0}^{N'} \#\{P \in \PDB \ | \ \varphi^n(P) \in \OS\},$$ so we're done by Corollary \ref{imagedensity}.
\end{proof}
\textbf{Acknowledgements:} It is a pleasure to thank Robert Grizzard and Aaron Levin for discussions related to the work in this paper.

\bibliographystyle{alpha}

\bibliography{bounded}
\end{document}